\documentclass{amsart}

\newtheorem{theorem}{Theorem}

\newtheorem{lemma}[theorem]{Lemma}

\newtheorem{conjecture}[theorem]{Conjecture}

\numberwithin{equation}{section}

\begin{document}

\title[The best possible quadratic refinement of Sendov's conjecture]
{The best possible quadratic refinement\\of Sendov's conjecture}
\author{Michael J. Miller}
\address{Department of Mathematics\\Le Moyne College\\Syracuse, New 
York 13214}
\email{millermj@lemoyne.edu}
\subjclass{Primary 30C15}
\keywords{Sendov, critical points, polynomial, derivative}
\thanks{21-Dec-2004} 

\begin{abstract}
A conjecture of Sendov states that if a polynomial has all its roots
in the unit disk and if $\beta$ is one of those roots, then within one
unit of $\beta$ lies a root of the polynomial's derivative.  If we
define $r(\beta)$ to be the greatest possible distance between $\beta$
and the closest root of the derivative, then Sendov's conjecture
claims that $r(\beta) \le 1$.

In this paper, we assume (without loss of generality) that $0 \le
\beta \le 1$ and make the stronger conjecture that $r(\beta) \le
1-(3/10)\beta(1-\beta)$.  We prove this new conjecture for all
polynomials of degree 2 or 3, for all real polynomials of degree 4,
and for all polynomials of any degree as long as all their roots 
lie on a line or $\beta$ is sufficiently close to $1$.
\end{abstract}

\maketitle

\section{Introduction}  

In 1958, Sendov conjectured that if a polynomial (with complex
coefficients) has all its roots in the unit disk, then within one unit
of each of its roots lies a root of its derivative.  A recent paper by
Sendov \cite{Sendov} and a recent book by Rahman and Schmeisser
\cite[Section 7.3]{Rahman-Schmeisser} both summarize the work that has
been done on this conjecture, identifying more than 80 related papers
that have been published in the past 35 years.  Despite this
substantial body of work, Sendov's conjecture has been verified only
for special cases.

Let $\beta$ be a complex number of modulus at most~$1$.  Define
$S(\beta)$ to be the set of polynomials of degree at least $2$ with
complex coefficients, all roots in the unit disk and at least one root
at $\beta$.  For a polynomial~$P$, define $d(P, \beta)$ to be the
distance between $\beta$ and the closest root of the derivative~$P'$.
Finally, define $r(\beta)=\sup \{ d(P, \beta): P \in S(\beta) \}$ and
note that $r(\beta) \le 2$ (since by the Gauss-Lucas Theorem
\cite[Theorem 2.1.1]{Rahman-Schmeisser} all roots of each $P'$ are
also in the unit disk, and so each $d(P, \beta)\le 2$).  In this
notation, Sendov's conjecture claims simply that $r(\beta) \le 1$.  To
date, the best such bound known to be true is that $r(\beta)\le
1.0753828$ \cite[Theorem~7.3.17]{Rahman-Schmeisser}.

In calculating $r(\beta)$, we will assume without loss of generality
(by rotation) that $0 \le \beta \le 1$.  Define $r_n(\beta)=\sup \{
d(P, \beta): P \in S(\beta) \hbox{ and } \deg P =n \}$.  Bojanov,
Rahman and Szynal have shown \cite[Lemma~4 and
$p(z)=z^n-z$]{Bojanov-Rahman-Szynal} that $r_n(0)=(1/n)^{1/(n-1)}$, so
letting $n$ tend to infinity gives $r(0) = 1$.  In addition,
Rubinstein has shown \cite[Theorem~1]{Rubinstein} that each
$r_n(1)=1$, so $r(1)=1$.  Given that $r(\beta)=1$ at both endpoints of
the interval $0 \le \beta \le 1$, the best possible linear (in
$\beta$) bound on $r(\beta)$ is that $r(\beta) \le 1$, which is the
claim of Sendov's conjecture.

To preserve the bounds of $1$ at $\beta=0$ and $\beta=1$, any
quadratic bound on $r(\beta)$ must be of the form $r(\beta) \le
1-c\beta(1-\beta)$ for some constant $c$.  Now recent work by Miller
\cite[Theorem 1, using $n=4$, $D_1=-1/5$ and $D_2=-2/5$]{Miller} shows
that $r_5(\beta)=1-(3/10)(1-\beta) +\mathcal O(1-\beta)^2$ for $\beta$
sufficiently close to $1$.  Matching slopes at $\beta=1$, we get $c=3/10$,
so the best possible quadratic bound on $r(\beta)$ is claimed by
\begin{conjecture}\label{3/10conj}
   For every $\beta \in [0,1]$ we have $r(\beta) \le
   1-(3/10)\beta(1-\beta)$.
\end{conjecture}

In addition to sharpening Sendov's conjecture, our Conjecture
\ref{3/10conj} has broader implications, as detailed in the following
three paragraphs.

Phelps and Rodriguez have conjectured \cite[after
Theorem~5]{Phelps-Rodriguez} that the only monic polynomials $P \in
S(\beta)$ with $d(P, \beta) \ge 1$ are of the form $z^n-e^{it}$ for
some~$t$.  Now Bojanov, Rahman and Szynal have shown \cite[Lemma
4]{Bojanov-Rahman-Szynal} that $d(P,0)<1$ for every $P \in S(0)$, so
if $P \in S(\beta)$ with $d(P, \beta) \ge 1$, then our
Conjecture~\ref{3/10conj} would imply that $\beta=1$.  Given this,
Rubinstein has shown \cite[Theorem 1]{Rubinstein} that
$P(z)=c(z^n-e^{it})$.  Thus our Conjecture~\ref{3/10conj} implies the
conjecture of Phelps and Rodriguez.

Since each $r_n(1)=1$, one might attempt to prove Sendov's conjecture
by showing that each $r_n(\beta)$ is increasing in the interval
$0<\beta<1$.  This cannot be done if our Conjecture~\ref{3/10conj} is
true, for then each $r_n(0.5) \le 0.925$, but for $n \ge 52$ we know
that $r_n(0)=(1/n)^{1/(n-1)}>0.925$.  Thus our
Conjecture~\ref{3/10conj} implies that for every $n \ge 52$, the bound
$r_n(\beta)$ must have a local minimum in the interval $0<\beta<1$.

Instead, one might try to show that $r_n(\beta)$ can never have a local
maximum in the interval $0<\beta<1$.  This also cannot be done in
general if our Conjecture~\ref{3/10conj} is true, for if we let
$C=\{
  1.98587$, $ 1.94369$, $ 1.87407$, $ 1.77800$, $ 1.65686$, 
$ 1.51240$, $ 1.34665$, $ 1.16197$, $ 0.96092$, $ 0.74623$, 
$ 0.52070$, $ 0.28718$, $ 0.04855$, $-0.19100$, $-0.42477$, 
$-0.64433$, $-0.84208$, $-1.01333$, $-1.15654$, $-1.27260$, 
$-1.36384$, $-1.43304$, \hfil\break
                        $-1.48297$, $-1.51587$, $-1.53341\}$, 
and define 
\begin{equation*}
P(z)=(z-0.09)(z+1)\prod_{c_i \in C}(z^2+c_i z+1),
\end{equation*}
then clearly $P \in S(0.09)$ and $\deg P=52$.  Calculating the
critical points of this polynomial (using at least 20 significant
digits) establishes that $r_{52}(0.09) \ge d(P, 0.09)>0.931$. Since
$r_{52}(0)=(1/52)^{1/51}\approx 0.925$, our Conjecture~\ref{3/10conj}
implies that $r_{52}(\beta)$ has a local maximum in the interval
$0<\beta<0.5$.  Furthermore, Borcea has conjectured
\cite[Remark~2.4]{Borcea} that each $r_n(\beta)$ is increasing in a
neighborhood of $\beta=0$.  If this is correct, then our
Conjecture~\ref{3/10conj} implies that $r_n(\beta)$ has a local
maximum in the interval $0<\beta<0.5$ for every $n \ge 52$.

\section{Results}  

In this paper, we will prove Conjecture \ref{3/10conj} for all
polynomials of degree 2 or 3, for all real polynomials of degree 4,
and for all polynomials of any degree as long as all their roots lie on
a line or $\beta$ is sufficiently close to $1$.  In most of these cases,
we will do even better than required for Conjecture \ref{3/10conj}, by
replacing $3/10$ with a larger number.

We begin by showing that Conjecture \ref{3/10conj} is true for all
polynomials of degree 2 or~3, via

\begin{theorem}\label{degree2&3}
   For every $\beta \in [0,1]$, we have $r_2(\beta) \le
   1-(1/2)\beta(1-\beta)$ and $r_3(\beta) \le 1-(1/3)\beta(1-\beta)$.
\end{theorem}
\begin{proof}
The result for $r_2(\beta)$ follows trivially from the fact that
$r_2(\beta)=(1+\beta)/2$. For polynomials of degree 3, Rahman has
shown \cite[Theorem 2]{Rahman} that
$r_3(\beta)=[3\beta+(12-3\beta^2)^{1/2}]/6$.  Since $0 \le 12-3\beta^2
\le [3+(1/2)(1-\beta^2)]^2$, then
\begin{equation*}
   r_3(\beta) \le [3\beta+3+(1/2)(1-\beta^2)]/6 \le
   1-(1/3)\beta(1-\beta),
\end{equation*}
and we are done.
\end{proof}

We next show that Conjecture \ref{3/10conj} is true provided that
$\beta$ is sufficiently close to~1 (where ``sufficiently close''
depends on the degree of the polynomial), using

\begin{theorem}\label{suffclose1}
   For every integer $n \ge 2$, if $\beta$ is sufficiently close to
   $1$, then we have $r_n(\beta) \le 1-(3/10)\beta(1-\beta)$.
\end{theorem}
\begin{proof}
Given Theorem \ref{degree2&3}, we may assume that $n \ge 4$.  Note that
\begin{equation*}
   1-\frac{3}{10}\beta(1-\beta)= 
   1-\frac{3}{10}(1-\beta)+\frac{3}{10}(1-\beta)^2.
\end{equation*}

If $n \ne 5$, then Miller has shown \cite[Theorem 1 and part 6 of
Lemma 8, using $n \ne 4$ and $c_{n+1}=D_1+D_2/n$]{Miller} that
$r_n(\beta) = 1+c_n(1-\beta)+\mathcal O(1-\beta)^2$ with $c_n <
-3/10$, so $r_n(\beta)<1-(3/10)(1-\beta)$ when $\beta$ is sufficiently
close to $1$ and the result follows.

If $n=5$, then Miller has shown \cite[Theorem 1, using $n=4$,
$D_1=-1/5$, $D_2=-2/5$, $D_3=0$, $D_4=0$, $D_5=1/25$ and
$D_6=-2/25$]{Miller} that 
\begin{equation*}
   r_5(\beta) = 1-\frac{3}{10}(1-\beta) +\frac{1}{200}(1-\beta)^2
   +\mathcal O(1-\beta)^3
\end{equation*}
and the result follows.
\end{proof}

We now examine polynomials having all real roots with
\begin{lemma} \label{allreal}
   If $P$ is a polynomial of degree $n \ge 2$ with all its roots in the 
   interval $[-1,1]$ and a root at $\beta$, then $d(P, \beta) \le 
   \max(2/n, 1/\sqrt{n})$.
\end{lemma}
\begin{proof}
Assuming (without loss of generality) that $P$ is monic, we can write
$P(z)=\prod_{i=1}^n (z-z_i)$ and $P'(z)=n\prod_{i=1}^{n-1}
(z-\zeta_i)$.  By Rolle's Theorem, there is a root of $P'$ between
every two roots of $P$, so we may order the roots so that
\begin{equation*}
   -1 \le z_1 \le \zeta_1 \le z_2 \le \zeta_2 \le \dots
                          \le z_{n-1} \le \zeta_{n-1} \le z_n \le 1.
\end{equation*}
For some $k$, we have $\beta=z_k$.  Then $P(z)=(z-\beta)\prod_{i=1,
i\ne k}^n (z-z_i)$, so $P'(\beta)=\prod_{i=1, i \ne k}^n (\beta-z_i)$
and so
\begin{equation} \label{equality}
   n\prod_{i=1}^{n-1} |\beta-\zeta_i|= |P'(\beta)| =
                           \prod_{i=1, i \ne k}^n |\beta-z_i|.
\end{equation}

If $k=1$, then $|\beta-\zeta_i| \ge |\beta-z_i|$ for $i=2, \dots,
n-1$.  Given equation \ref{equality}, this implies that
$n|\beta-\zeta_1| \le |\beta-z_n|\le2$, and so $d(P,\beta)
= |\beta-\zeta_1| \le 2/n$.

If $1<k<n$, then $|\beta-\zeta_i|\ge|\beta-z_{i+1}|$ for $i=1,
\dots, k-2$ and that $|\beta-\zeta_i| \ge |\beta-z_i|$ for $i=k+1, \dots,
n-1$.   Given equation \ref{equality}, this implies that
$n|\beta-\zeta_{k-1}||\beta-\zeta_k| \le
|\beta-z_1||\beta-z_n| \le (1-\beta)(1+\beta) \le 1$, so $d(P, \beta)
= \min(|\beta-\zeta_{k-1}|,|\beta-\zeta_k| ) \le 1/\sqrt{n}$.

If $k=n$, then $|\beta-\zeta_i| \ge |\beta-z_{i+1}|$ for $i=1, 
\dots, n-2$.   Given equation \ref{equality}, this implies that
$n|\beta-\zeta_{n-1}| \le |\beta-z_1| \le 2$, and so $d(P, \beta) =
|\beta-\zeta_{n-1}| \le 2/n$.
\end{proof}

We now show that Conjecture \ref{3/10conj} is true for polynomials
with all roots on a line, via
\begin{theorem} \label{rootsonline}
   If all the roots of $P \in S(\beta)$ lie on a line, then 
 $d(P, \beta) \le 1-(1/2)\beta({1-\beta})$.
\end{theorem}
\begin{proof}
Let $n=\deg P$.  Given Theorem \ref{degree2&3}, we may assume that $n
\ge 3$.  Moving the roots of $P$ to the interval $[-1,1]$ by a rigid
transformation of the plane (and thus leaving $d(P,\beta)$ unchanged),
we see by Lemma \ref{allreal} that $d(P, \beta) \le \max(2/n,
1/\sqrt{n}) \le 2/3$.  Our result follows from the observation that
$1-(1/2)\beta(1-\beta)$ has a minimum of $7/8$ for $0 \le \beta \le
1$.
\end{proof}

To examine real polynomials of degree $4$, we will need
\begin{lemma}\label{modpprime}
   For every monic real polynomial $P \in S(\beta)$ of degree $4$ with
   $d(P,\beta)>(1+\beta)/2$, we have $|P'(\beta)| \le (1+\beta)^2$.
\end{lemma}
\begin{proof}
Write $P(z)=(z-\beta)\prod_{i=1}^3 (z-z_i)$ and note that
$P'(\beta)=\prod_{i=1}^3 (\beta-z_i)$.

If $P$ has a root (say) $z_1$ in the half-plane $\{z: \Re(z) \ge \beta
\}$, then $|z_1-\beta| \le 1$.  For $i \ge 2$, we have $|z_i-\beta| \le
1+\beta$, and so $|P'(\beta)|=\prod_{i=1}^3 |\beta-z_i|\le(1+\beta)^2$
and we are done.

Assume then that the half-plane $\{z: \Re(z) \le \beta \}$ contains
all roots of $P$ and hence also (by the Gauss-Lucas theorem
\cite[Theorem 2.1.1]{Rahman-Schmeisser}) all roots of $P'$. Since by
hypothesis $d(P,\beta)>(1+\beta)/2$, then any real roots of $P'$ would
be in the interval $[-1, (\beta-1)/2)$.  Now $P$ is a real polynomial
of even degree with a real root $\beta$, hence $P$ has another real
root (say) $z_1 \ge -1$.  Since $P'$ has no real roots in the interval
$[(\beta+z_1)/2,\infty)$ and since $P'(x)>0$ for large $x$, then
$P'((\beta+z_1)/2)>0$.  Now
\begin{equation*}
   P'(z)=(z-\beta)(z-z_1)(2z-z_2-z_3)+(2z-\beta-z_1)(z-z_2)(z-z_3),
\end{equation*}
so 
\begin{equation*}
   0 < P'\left(\frac{\beta+z_1}{2}\right)
   = -\frac{(\beta-z_1)^2}{4}(\beta+z_1-z_2-z_3),
\end{equation*}
and so we have $z_2+z_3 > \beta+z_1 \ge \beta-1$.  Now $z_2$ and $z_3$
cannot both be real, else the interval $((\beta-1)/2,\beta]$ would
contain the larger of the two as well as $\beta$, and hence by Rolle's
theorem a root of $P'$.  Thus $z_2$ and $z_3$ must be complex
conjugates with $\Re(z_2)=\Re(z_3)>(\beta-1)/2$, so for $i \ge 2$ we
have
\begin{equation*}
   |z_i-\beta|^2 = z_i \bar z_i - 2\beta\Re(z_i) + \beta^2 
   \le 1 -2\beta(\beta-1)/2 + \beta^2 = 1+\beta.
\end{equation*}
Note that $|z_1-\beta| \le 1+\beta$, and so $|P'(\beta)|=\prod_{i=1}^3
|\beta-z_i|\le(1+\beta)^2$ and we are done.
\end{proof}

We mention in passing that Lemma \ref{modpprime} may fail for nonreal
polynomials, as can be seen by choosing $\beta=0.674$ and
$P(z)=\int_\beta^z 4(w+0.24-0.38i)(w+0.13+0.25i)^2\,dw$.  A numerical
calculation establishes that the roots of $P$ have moduli less than
$1$, so $P \in S(\beta)$.  However, $d(P,\beta) \approx 0.84197 >
(1+\beta)/2$, but $|P'(\beta)| \approx 2.80687 > (1+\beta)^2$.

Finally, we prove Conjecture \ref{3/10conj} for all real polynomials
of degree 4 with
\begin{theorem}\label{degree4real}
   For every real polynomial $P \in S(\beta)$ of degree $4$ we have
   $d(P, \beta) \le 1-(1/3)\beta(1-\beta)$.
\end{theorem}
\begin{proof}
If $d(P, \beta) \le (1+\beta)/2$, then $d(P, \beta) \le
1-(1/2)\beta(1-\beta)$ and we are done.  Assume then that $d(P,
\beta)>(1+\beta)/2$, assume (without loss of generality) that $P$ is
monic and write $P'(z)=4\prod_{i=1}^3 (z-\zeta_i)$.  Then using Lemma
\ref{modpprime} we have
\begin{equation}\label{prod}
   4(d(P, \beta))^3 \le 4\prod_{i=1}^3 |\zeta_i-\beta| = |P'(\beta)|
   \le (1+\beta)^2.
\end{equation}

Expanding the cube shows that $1-x \le (1-x/3)^3$ for $0 \le x \le 1$,
so letting $x=\beta(1-\beta)$ gives us
\begin{equation}\label{cube}
   (1+\beta)^2 \le 4(1-\beta(1-\beta)) \le 4(1-\beta(1-\beta)/3)^3.
\end{equation}

Combining lines \ref{prod} and \ref{cube} gives us that $4(d(P,
\beta))^3 \le 4(1-\beta(1-\beta)/3)^3$ and our result follows.
\end{proof}


\end{document}